\documentclass[12pt]{amsart}
\usepackage{amsfonts,latexsym,amsthm}

\input epsf

\newtheorem{theorem}{Theorem}[section]
\newtheorem{lemma}[theorem]{Lemma}
\newtheorem{corol}[theorem]{Corollary}
\newtheorem{prop}[theorem]{Proposition}

{\theoremstyle{definition} \newtheorem{defin}{Definition}[section]}
{\theoremstyle{remark} 
\newtheorem{example}{Example}[section]}

\newcommand{\Cbb}{{\mathbb{C}}}

\newcommand{\Pbb}{{\mathbb{P}}}

\newcommand{\Qbb}{{\mathbb{Q}}}
\newcommand{\Til}[1]{{\widetilde{#1}}}

\newcommand\csm{{c_{\text{SM}}}}
\newcommand\cma{{c_{\text{Ma}}}}
\newcommand\sma{{s_{\text{Ma}}}}
\newcommand\cfu{{c_{\text{F}}}}
\newcommand\cwm{{c_{\text{wMa}}}}
\DeclareMathOperator{\PGL}{PGL}

\title{Interpolation of characteristic classes of singular hypersurfaces}
\author{Paolo Aluffi and Jean-Paul Brasselet}
\date{July 2001}
\address{Mathematics Department, Florida State University,
Tallahassee FL 32306, U.S.A.}
\email{aluffi@math.fsu.edu}
\address{
IML -- CNRS, 
Case 907, Luminy, 
13288 Marseille Cedex 9, France
}
\email{jpb@iml.univ-mrs.fr}

\begin{document}

\begin{abstract}
We show that the Chern-Schwartz-MacPherson class of a hypersurface
$X$ in a nonsingular variety $M$ `interpolates' between two other
notions of characteristic classes for singular varieties, provided
that the singular locus of $X$ is smooth and that certain numerical invariants 
of $X$ are constant along this locus. This allows us to define a lift 
of the Chern-Schwartz-MacPherson class of such `nice' hypersurfaces to 
intersection homology. As another application, the interpolation result 
leads to an explicit formula for the Chern-Schwartz-MacPherson class of $X$ 
in terms of its polar classes.
\end{abstract}

\subjclass
{Primary 14C17;
Secondary 57R20 32C40}

\maketitle

\section{Introduction and statement of the result}\label{intro}

There are several different notions of `characteristic classes' of
possibly singular varieties, generalizing the notion of (homology)
Chern classes of the tangent bundle of nonsingular ones; the
relationship between some of these classes has been the object of
recent work. In this note we prove a formula relating the
Chern-Mather class, the Chern-Schwartz-MacPherson class, and the
class of the virtual tangent bundle of a hypersurface in a nonsingular
variety, under the assumption that the singular locus of $X$ is smooth
and that certain numerical invariants of $X$ are constant along this
locus. As an application, we obtain an explicit formula for
Chern-Schwartz-MacPherson's class of $X$ in terms of its polar
classes (generalizing a result in \cite{5}), under the same hypothesis
on its singularity, and assuming that $X$ is quasi-projective.

A more immediate, but perhaps more striking, application is to the problem 
of lifting Chern-Schwartz-MacPherson's classes of a hypersurface to 
intersection homology. While the class of the virtual tangent bundle of $X$ 
trivially has a natural lift to intersection homology, examples of Mark Goresky 
and Jean-Louis Verdier show that the problem of lifting 
Chern-Schwartz-MacPherson's classes is much subtler (see \cite{5.3} for a 
discussion of these examples). A lift exists to intersection homology 
with rational coefficients in middle perversity as a consequence of 
\cite{3.5}; but there is no known way to construct a `canonical' such lift in 
general. For quasi-projective hypersurfaces $X$ satisfying our hypothesis, the 
interpolation formula given below defines a lift of the 
Chern-Schwartz-MacPherson class of $X$ in $IH_*(X)$ with rational 
coefficients, in middle perversity. This could be used to define {\em Chern 
numbers\/} $c_i\,c_{\dim X-i}$ for projective singular hypersurfaces satisfying 
the condition considered here. Computing such numbers explicitly would be very 
interesting; also it would be interesting to establish the exact dependence 
(if any) of these numbers or of our lift on the embedding of $X$.

Let $X$ be a reduced hypersurface of a nonsingular complex algebraic
variety $M$. We denote by
$$\csm(X)\quad,\quad \cma(X)\quad,\quad \cfu(X)$$
the three classes mentioned above; the first two are defined in
\cite{9}, while the third is the class of the virtual tangent bundle of 
$X$:
$$\cfu(X)=c\left(\frac{TM}{\mathcal O(X)}\right)\cap [X]=c(TM)\cap
s(X,M)\quad.$$
Here, $s(X,M)$ is the {\em Segre class\/} of $X$ in
$M$, cf.~\cite{6}, Chapter~4. The class $\cfu(X)$ equals William
Fulton's intrinsic class of $X$, which can be defined for {\em every\/}
scheme embeddable in a nonsingular variety (cf.~\cite{6}, \S4.2.6), and is 
independent of the ambient variety $M$.

All these classes can be defined in a good homology theory on $X$; in
this paper we work in the Chow group of $X$ with rational coefficient,
denoted $(A_*X)_{\Qbb}$ (except for the application to lifting to $IH_*$).

We denote by $Y$ the singularity {\em subscheme\/} of $X$ (locally defined 
by the partial derivatives of a generator of its ideal in $M$) and
by $Y'$ its support, that is, the singular {\em locus\/} of $X$. For $p\in Y'$ 
we consider two numerical invariants of $X$ at $p$: 

---the local Euler obstruction of $X$ at $p$, $\text{Eu}_X(p)$, and

---the Euler characteristic $\chi_p$ of the Milnor fiber of $X$ at
$p$.

\begin{defin}\label{nice}
A variety $X$ is a {\em nice hypersurface\/} if it can be realized as a 
hypersurface in a nonsingular variety $M$, and further its singular locus 
$Y'$ is nonsingular and irreducible, and the numbers $\text{Eu}_X(p)$, 
$\chi_p$ are constant along $Y'$.
\end{defin}

{\em In the results given below, we assume that $X$ is a nice 
hypersurface,\/} and we denote the constant values of $\text{Eu}_X(p)$, 
$\chi_p$ by $\text{Eu}$, $\chi$ respectively. This condition is satisfied 
for example if the stratification $\{Y',X\setminus Y',M\setminus X\}$ of 
the ambient variety $M$ is Whitney regular, or satisfies the weaker condition 
of {\em $c$-regularity\/} of Karim Bekka, \cite{4}\footnote{We are indebted 
to J\"org Sch\"urmann for this remark; see for example \cite{13}, V.3 and 
VI.2.}. The precise algebro-geometric requirement is that the normal cone of the 
singularity subscheme $Y$ in $M$ be irreducible, cf.~Lemma~\ref{lemma3}.

Under this assumption, we will show that the Chern-Schwartz-MacPherson
class of $X$ is an `interpolation' of the classes $\cma(X)$ and
$\cfu(X)$; we will now state this result precisely. As $X$ is a 
hypersurface in $M$, it determines a line bundle $\mathcal O_M(X)$; we adopt 
a common abuse of notation and denote by $X$ the first Chern class $c_1(\mathcal 
O_M(X))$ of this line bundle, and its restrictions to subschemes of $M$. Thus 
$s(X,M)=\frac{[X]}{1+X}$ in our notations, and $\frac 1{1+\alpha X}$ is 
shorthand for $1-\alpha\, c_1(\mathcal O(X))+\alpha^2\, c_1(\mathcal O(X))^2+
\cdots$.

For all rational numbers $\alpha$, we let
$$c_{(\alpha)}(X)=\cfu(X)+\frac{(1-\alpha)}{1+\alpha X}\left(\cma(X) 
-\cfu(X)\right) \quad\in(A_*X)_{\Qbb}\quad.$$
Thus 
$$c_{(0)}(X)=\cma(X)\quad\text{and}\quad c_{(1)}(X)=\cfu(X)$$
trivially. Also, $c_{(\alpha)}(X)$ does {\em not\/} depend on the ambient 
manifold $M$ in which $X$ is realized as a hypersurface: indeed, the class 
$(\cma(X) -\cfu(X))$ is supported on the singular locus $Y'$ of $X$, and it 
can be shown that the action of $X=c_1(\mathcal O_M(X))$ on $Y'$ is 
independent of $M$. In fact, the restriction of $\mathcal O_M(X)$ to $Y$ 
does not depend on $M$.
\begin{theorem}\label{main}
If $X$ is a nice hypersurface, then
$$c_{(\rho)}(X)=\csm(X)\quad\text{in $(A_*X)_{\Qbb}$,}$$
where $\csm(X)$ denotes the Chern-Schwartz-MacPherson class of $X$, and
$\rho=\frac{1-\text{Eu}}{\chi-\text{Eu}}$.
\end{theorem}

We note that, under our assumption on the singularity of $X$, and by the 
very definition of $\csm(X)$ in \cite{9}, the class $\csm(X)$ is a simple 
linear combination of $\cma(X)$ and of the total homology Chern class of the 
singular locus $Y'$ of $X$. Our formula replaces this latter {\em local\/} 
ingredient with the {\em global\/} information of the class of the virtual 
tangent bundle of $X$. If the Chern classes of $Y'$ are known, our formula can 
be used `in reverse' to obtain information about the invariants Eu, $\chi$, 
bypassing a local study of $X$ near its singular locus. See \S\ref{remex}, for 
a precise statement (Proposition~\ref{proposition7}) and an example of such 
a computation.

Also note that while the condition of being `nice' is of course very 
strong in general, it is automatically satisfied in codimension equal to 
the codimension of the singular locus $Y'$ of $X$, provided that the 
top-dimensional part of $Y'$ is irreducible: indeed, the hypersurface obtained 
from $X$ by removing the locus where the invariants jump is then trivially 
nice. The formula of the theorem gives then a relation between the 
classes up to that codimension, under the sole assumption that the 
singular locus is irreducible in top dimension; see Example~\ref{example4.2} 
for an illustration of this fact.

As mentioned earlier, the interpolation formula can be used to lift 
Chern-Schwartz-MacPherson's classes of a quasi-projective nice 
hypersurface $X$ to $IH_*(X)$ (with rational coefficients, in middle 
perversity). This relies precisely on our trading the `local' information 
of $c(TY')\cap [Y']$ (which is hard to transfer into $IH_*(X)$) for the 
`global' information of $\cfu(X)$. In order to define the lift, we just 
note that a lift of $\cma(X)$ for quasi-projective $X$ is defined in 
\cite{5.7}; $\cfu(X)$ lifts as it is the class of the virtual tangent bundle 
of $X$; and the other ingredients in the interpolation formula also involve 
elements in cohomology, so they trivially lift to $IH_*(X)$. Hence the formula 
defines an element of $IH_*(X)$ for nice singular hypersurfaces, 
which lifts $\csm(X)$ by the main theorem. We note that as the lift of 
$\cma(X)$ of \cite{5.7} potentially depends on the realization of $X$ as a 
quasi-projective variety, our lift of $\csm(X)$ may also depend on this 
choice. It would be interesting to establish whether it is in fact 
uniquely determined by $X$ itself.

Theorem~\ref{main} is proved in \S\ref{segre}; the main tools are formulas 
from \cite{2}, manipulations of Segre classes, and a key result of Adam 
Parusi{\'n}ski and Piotr Pragacz (\cite{12}). In \S\ref{relpol} we consider a 
situation in which $\cma(X)$ can be expressed very concretely, that is, 
when $X$ is given explicitly as a subvariety of a projective space $\Pbb^n$. 
In this case, a result of Ragni Piene can be used to express the class in
terms of a suitable combination $[P]$ of the {\em polar classes\/} of $X$:
$$[P]=-\sum_{k\ge 0} [P_k]^\vee\otimes_{\Pbb^n}\mathcal O(1)\quad,$$
where $[P_k]$ denotes the class of the $k$-th polar locus of $X$, and
we use the notations introduced in \cite{1}. We recall precise
definitions and Piene's result in \S\ref{relpol}. We then have:
\begin{corol}\label{polar}
If $X\subset \Pbb^n$, and $X$ is a nice
hypersurface (in some variety $M$), let
$$\rho=\frac{1-\text{Eu}}{\chi-\text{Eu}}\quad,\quad 
\sigma=1-\rho=\frac{\chi-1}{\chi-\text{Eu}}\quad.$$
Then the Chern-Schwartz-MacPherson class of $X$ is given by
$$\csm(X)=c(TM)\cap \frac{\rho\,[X]}{1+\rho X}+c(T\Pbb^n)\cap
\frac{\sigma\, [P]}{1+\rho X}\in (A_*X)_{\Qbb}\quad.$$
\end{corol}

As in Theorem~\ref{main}, $X$ is used here to denote $c_1(N_X M)$.
Note that we are not requiring $X$ to be a hypersurface {\em in
$\Pbb^n$;\/} all we need is that $X$ can be abstractly
realized as a nice hypersurface in {\em some\/} variety $M$, and that $X$ 
is itself quasi-projective.

If $Y=Y'$ and the (nice) hypersurface $X$ has multiplicity~2 along $Y$, 
then $(\chi,\text{Eu})=(2,0)$ if the codimension of $Y$ in $X$ is even, 
and $(0,2)$ if it is odd (this follows from Lemma~\ref{lemma3} in 
\S\ref{segre}). In both cases $\rho=\sigma=\frac 12$; if further $M=\Pbb^n$, 
then the formula in the corollary specializes to the case considered in 
\cite{5}; this was our starting point in this work.

Other expressions for Chern-Schwartz-MacPherson's class of a singular 
variety $X$ in the context of the study of polar varieties are known, 
notably those given in \S6 of \cite{8} (without any restriction on $X$!). 
The work of Ragni Piene, L\^e D\~ung Tr\'ang, and Bernard Teissier (\cite{8},
\cite{11}) has exposed the close relationship between characteristic classes 
of singular varieties and their polar varieties.
Our motivation in \S3 is however somewhat different than in these
references---we have aimed specifically at identifying the
contribution of polar varieties to `correction terms' between
different notions of characteristic classes. The term $c(TM)\cap
\frac{\rho[X]}{1+\rho X}$ is the analog of Fulton's intrinsic class
for a `virtual' hypersurface $\rho X$. The other term, $c(T\Pbb^n)\cap
\frac{\sigma[P]}{1+\rho X}$, could then be interpreted as a {\em
Milnor class\/} (in the sense of \cite{5}) for such a virtual
hypersurface. Thus, Corollary~\ref{polar} brings evidence to
the possibility that Milnor classes admit simple expressions in terms
of polar classes. Positive results in this direction could lead to a
good treatment of Milnor classes for more general varieties, which
would be highly desirable.

{\bf Acknowledgements.\/} We are grateful to Ragni Piene and J\"org 
Sch\"urmann for useful comments. The first author would like to thank the 
IML, Marseille and the French CNRS for the generous hospitality while this 
work was completed.

\section{Segre classes and characteristic classes}\label{segre}

Let $X$ be a reduced hypersurface of an arbitrary nonsingular variety
$M$; we work over $\Cbb$ for convenience, but most of what we say
can easily be extended to arbitrary algebraically closed fields of
characteristic~0. 

We denote by $Y$ the singularity subscheme of $X$, and for the
moment we make no further assumptions on $X$ or $Y$. The following
formulas for the Mather and {Schwartz-MacPher\-son} classes (denoted
respectively $\cma(X)$, $\csm(X)$) are given in \cite{2} (Lemma~I.2,
Theorem~I.3):
\begin{lemma}\label{lemma1}
Let $\pi:\Til M \to M$ be the blow-up of $M$ along
$Y$; let $\mathcal Y$ be the exceptional divisor in this blow-up, and let
$\mathcal X$, $\Til X$ respectively be the total and strict transforms of $X$
in $\Til M$. Then
$$\aligned
\cma(X) &=c(TM)\cap \pi_*\left(\frac{[\Til X]}{1+\mathcal X-\mathcal Y}\right)\\
\csm(X) &=c(TM)\cap \pi_*\left(\frac{[\mathcal X-\mathcal Y]}{1+\mathcal X-
\mathcal Y}\right)\quad.
\endaligned$$
\end{lemma}
These formulas can be conveniently rewritten without reference to
classes in $\Til M$, by adopting the notations introduced in \cite{1}:
for $A=\sum a_p$, $a_p$ a class of dimension $p$ in a subscheme of
$M$, and $\mathcal L$ a line bundle, we let $A \otimes \mathcal L$ be the
class $\sum_p c(\mathcal L)^{p-\dim M} \cap a_p$ and $A^\vee$ be $\sum
(-1)^{p-\dim M}a_p$. Note that these notations depend on the ambient
variety $M$; this will be understood in the following. With these
notations:
\begin{prop}\label{prop2}
 Denote by $\mathcal L$ the restriction of the line
bundle $\mathcal O(X)$ to $Y$. Then
$$\aligned
\cma(X) &=c(TM)\cap \left(\frac{[X]}{1+X}+s(Y,X)^\vee\otimes\mathcal
L\right)\\
\csm(X) &=c(TM)\cap \left(\frac{[X]}{1+X}+(c(\mathcal L)\cap
s(Y,M))^\vee\otimes\mathcal L\right)
\endaligned$$
\end{prop}
\begin{proof}
The statement follows at once from the
formulas given in Lemma~\ref{lemma1}, with standard manipulations involving the
notations recalled above. For the first formula, write
$$\frac{[\Til X]}{1+\mathcal X-\mathcal Y}=\frac{[\Til X]}{1-\mathcal Y}\otimes 
\mathcal L
=\left([\Til X]+\frac{\mathcal Y\cdot[\Til X]}{1-\mathcal Y}\right)\otimes
\mathcal
L=\frac{[\Til X]}{1+\mathcal X}+\left(\frac{\mathcal Y\cdot[\Til X]}{1+
\mathcal Y}
\right)^\vee\otimes\mathcal L\quad:$$
pushing this expression forward by $\pi$ gives
$$\frac{[X]}{1+X}+s(Y,X)^\vee\otimes\mathcal L\quad,$$
yielding the first formula in the statement. The second formula is
proven similarly; it is in fact Theorem~I.4 in \cite{2}.
\end{proof}
Proposition~\ref{prop2} highlights the difference between the two notions 
of Chern-Mather and Chern-Schwartz-MacPherson classes of a hypersurface $X$ 
in a nonsingular variety~$M$: the distinction lies in the difference
between
$$s(Y,X)\quad\text{and}\quad c(\mathcal O(X))\cap s(Y,M)$$
and in this sense it is precisely captured by the singularity
subscheme $Y$ of the hypersurface. Relating the two characteristic
classes directly amounts then to comparing $s(Y,X)$ and $s(Y,M)$
directly. Unfortunately, very few such comparison results are known in
general (cf.~\cite{6}, Example~4.2.8 for a counterexample to the
naive guess for such a comparison). At present, the strong assumption
posed in \S\ref{intro} to state the main theorem of this paper is necessary
precisely because it allows us to perform this comparison.

First we gather more information from the blow-up $\Til M$; here we make 
crucial use of a result from \cite{12}.
\begin{lemma}\label{lemma3}
Under the hypotheses of the theorem, $\mathcal Y$ is
irreducible. Denoting by $\mathcal Y'$ its support, and writing
$$\mathcal Y =m\mathcal Y'\quad, \quad \mathcal X =\Til X+n\mathcal Y'$$
as cycles, then we have
$$m=(-1)^{\dim X-\dim Y}(\chi-1)\quad, \quad n=(-1)^{\dim X-\dim
Y}(\chi-\text{Eu})\quad.$$
with $\chi$ and $\text{Eu}$ defined as in \S\ref{intro}.
\end{lemma}
\begin{proof} 
Let $\mathcal Y=\sum m_i\mathcal Y_i$ be the irreducible
decomposition of $\mathcal Y$; by \cite{7}, each $\mathcal Y_i$ can
be identified with the conormal variety of its support. In particular,
there is exactly one component $\mathcal Y'$ over the support $Y'$ of $Y$.
Following \cite{12} we let $\mu=(-1)^{\dim X}(\chi -1)$, and remark
that under our assumption this is a multiple of the characteristic
function $1_Y$ of $Y'$, hence of the local Euler obstruction
$\text{Eu}_{Y'}$ since $Y'$ is nonsingular by hypothesis. By
Theorem~2.3(iii) of \cite{12}, the cycle $\sum m_i\mathcal Y_i$ must
equal a constant times $\mathcal Y'$; it follows that $\mathcal Y'$ is the
only irreducible component of $\mathcal Y$, and further that 
$$\mu=m (-1)^{\dim Y}1_Y\quad.$$ 
The first assertion in the statement follows, as well as the formula
for $m$. The formula for $n$ can be obtained similarly from
Theorem~2.3 in \cite{12}, or from the fact that
$$1_X=\text{Eu}_X(p)+(n-m) (-1)^{\dim X-\dim Y}\text{Eu}_Y(p)$$
for all $p\in Y$ (as proved in \cite{3}, \S2).
\end{proof}

Note that, for nice hypersurfaces, necessarily $\chi\ne 1$ and $\chi\ne 
\text{Eu}$ (this follows from Lemma~\ref{lemma3}). The next lemma relates 
$s(Y,X)$ and $s(Y,M)$ in the special case of nice hypersurfaces.

\begin{lemma}\label{lemma4}
 With assumptions and notations as above,
$$s(Y,X)=\left(\frac {\chi-\text{Eu}}{\chi-1}+X\right)\cdot s(Y,M)$$
in $(A_*Y)_{\Qbb}$.
\end{lemma}
\begin{proof}
Since $Y$ has codimension at least~2 in $M$, note that
$$\pi_*\frac{\mathcal Y\cdot [\mathcal Y]}{1+\mathcal Y}=\pi_*\frac{[\mathcal
Y]}{1+\mathcal Y}=s(Y,M)\quad.$$ 
Therefore with notations as in Lemma~\ref{lemma3}
$$\aligned
s(Y,X)&=\pi_*\frac{\mathcal Y\cdot[\Til X]}{1+\mathcal Y}= \pi_*\frac{\mathcal Y
\cdot([\mathcal X]-n[\mathcal Y'])}{1+\mathcal Y}=X\cdot\pi_*\frac{[\mathcal Y]}{1+\mathcal Y}
- \frac nm\pi_*\frac{\mathcal Y\cdot[\mathcal Y]}{1+\mathcal Y}\\
&=X\cdot s(Y,M)+\frac nm s(Y,M)=\left(\frac nm+X\right)\cdot s(Y,M) 
\endaligned$$
\end{proof}

After these preliminary considerations, we are ready to prove the main
theorem. Mimicking the relation between $s(X,M)$ and $\cfu(X)$, we write 
$\sma(X,M)$ for $c(TM)^{-1}\cap \cma(X)$; and, as above, 
$\sigma=\frac{\chi-1}{\chi-\text{Eu}}$ and $\rho=1-\sigma$.
\begin{proof}[Proof of Theorem~\ref{main}]
In view of the first formula in
Proposition~\ref{prop2}, 
$$\sma(X,M)=\frac{[X]}{1+X}+s(Y,X)^\vee\otimes\mathcal L\quad.$$
By Lemma~\ref{lemma4}, $s(Y,X)=\frac 1\sigma (1+\sigma X)\cap s(Y,M)$,
hence
$$s(Y,M)=\sigma\frac{1}{1+\sigma X}\cap s(Y,X)\quad.$$
Therefore
$$\aligned
\left(c(\mathcal L)\cap s(Y,M)\right)^\vee \otimes\mathcal L &= \sigma\left(
\frac{1-X}{1-\sigma X}\cap s(Y,X)^\vee \right)\otimes \mathcal L\\
& = (1-\rho) \frac 1{1+\rho X}\cap (s(Y,X)^\vee \otimes \mathcal L)
\endaligned$$
(using Proposition~1 in \cite{1}). That is,
$$\left(c(\mathcal L)\cap s(Y,M)\right)^\vee \otimes\mathcal L =
\frac{(1-\rho)}{1+\rho X}\left(\sma(X,M)-\frac{[X]}{(1+X)}\right)
\quad,$$
and the formula given in the theorem now follows from the second formula
of Proposition~\ref{prop2} and from $\cfu(X)=c(TM)\cap\frac{[X]}{1+X}$.
\end{proof}

\section{Relation with polar classes}\label{relpol}

We now turn to the situation in which $X$ is embedded in a projective
space $\Pbb^n$, and to polar classes; all we need to do is rewrite a
result of Ragni Piene into our language. If $X$ is a (closed) subvariety 
of dimension $r$ in $\Pbb^n$, the $k$-th {\em polar locus\/} $P_k$ of $X$
with respect to a general linear subspace $L_k= \Pbb^{k-2+n-r}\subset
\Pbb^n$ is the closure of the locus
$$\{x\in X_{\text{smooth}}\, |\, \dim(T_x X\cap L_k)\ge k-1\}\quad.$$
The class $[P_k]\in A_{\dim X-k}$ is independent of the (general)
choice of the subspace $L_k$ (cf.~\cite{10}, Proposition~1.2).
Note that $[P_0]=[X]$, since the condition defining $P_0$ is vacuous.
We define the `total polar class' of $X$ by
$$[P]=(-1)^{n-r}\sum_{k\ge 0} [P_k]^\vee\otimes_{\Pbb^n}\mathcal O(1)\quad;$$
again we are using the notations of \cite{1}, for convenience.

The main observation here is that the class $[P]$ is closely related to the
class $\cma(X)$. The precise relation is given in part (a) in the following 
theorem, due to Ragni Piene; we include part (b) for completeness, and stress 
that (b) holds for {\em arbitrary\/} hypersurfaces.
\begin{theorem}\label{theorem5}
(a) (R.~Piene) For any subvariety $X$ of $\Pbb^n$
as above:
$$\cma(X)=c(T\Pbb^n)\cap [P]\quad.$$
(b) If, further, $X$ is a hypersurface in a nonsingular variety $M$,
with singularity subscheme $Y$, and $\mathcal L$ denotes $\mathcal
O_M(X)_{|X}$, then
$$s(Y,X)=[X]+\frac{c(N_X^*\Pbb^n\otimes\mathcal L)}{c(\mathcal L)^{n-r-1}}
\cap\left([P]^\vee\otimes_M\mathcal L\right)\quad\in A_*X
\quad.$$
\end{theorem}

\begin{proof} 
(a) This is the translation in our notations of the second formula 
in Ragni Piene's Th\'eor\`eme~3 in \cite{11}:
$$\cma(X)=\sum_{k\ge 0} \sum_{i=0}^k (-1)^{k-i} \binom{r+1-k+i}{i}
H^i\cdot [P_{k-i}]\quad,$$
where $H$ is the hyperplane class.

(b) From part (a) and Proposition~2, we have
$$c(T\Pbb^n)\cap [P]=c(TM)\cap \left(\frac{[X]}{1+X}+s(Y,X)^\vee
\otimes\mathcal L\right)\quad;$$
therefore, using Propositions~1 and 2 from \cite{1}:
$$\aligned
\frac{[X]}{1+X}+s(Y,X)^\vee\otimes \mathcal L &=\frac{c(T\Pbb^n)}{c(TM)}
\cap[P] =\frac{c(N_X\Pbb^n)}{c(\mathcal L)}\cap [P]\\
s(Y,X)^\vee\otimes \mathcal L &=\frac 1{c(\mathcal L)}\left(c(N_X\Pbb^n)\cap
[P]-[X]\right)\\
s(Y,X)\otimes \mathcal L^\vee &=\frac 1{c(\mathcal L^\vee)}
\left(c(N^*_X\Pbb^n)\cap [P]^\vee+[X]\right)\\
s(Y,X) &=c(\mathcal L)\frac {c(N^*_X\Pbb^n\otimes\mathcal L)}
{c(\mathcal L)^{n-r}}\cap\left([P]^\vee\otimes\mathcal L\right)+[X]\quad,
\endaligned$$
with the stated result.
\end{proof}
In the particular case in which $X$ is a hypersurface {\em of
$\Pbb^n$\/}, the formula in part (b) reduces to
$$s(Y,X)=[P]^\vee\otimes\mathcal L+[X]\quad;$$
as the reader may check, this is equivalent to Piene's {\em Pl\"ucker
formulae\/} (cf.~Theorem~2.3 in \cite{10}).

Corollary~\ref{polar} follows from part (a), Theorem~\ref{main}, and 
straightforward manipulations.

\section{Remarks and examples}\label{remex}

Regarding the computability of the key coefficient $\rho$ needed in
order to apply Theorem~\ref{main}, the following observation may be useful.
\begin{prop}\label{proposition7}
With notations and assumptions as in \S~\ref{segre},
$$(1+X)\left(\cma(X)-\cfu(X)\right)=\left((\text{Eu}-\chi) + 
(\text{Eu}-1) X\right)\cdot (c(TY')\cap [Y'])\quad.$$
\end{prop}
\noindent The point is that if $\cma(X)$, $\cfu(X)$, and $c(TY')$ are
known, and $X\cdot\dim Y'\ne 0$ (for example, $\dim Y'>0$ if
$M=\Pbb^n$), then this formula determines $(\text{Eu}-1)$
and $(\text{Eu}-\chi)$; $\rho$ is the quotient of these two numbers.
\begin{proof}
Since $Y'$ is assumed to be nonsingular, the weighted
Chern-Mather class of $Y$ (cf.~\cite{3} and Lemma~\ref{lemma3}) is given by
$$\cwm(Y)=m c(TY')\cap[Y']=(-1)^{\dim X-\dim Y}(\chi-1)
c(TM)\cap s(Y',M)\quad;$$
on the other hand, by Proposition~1.3 in \cite{3}
$$\aligned
\cwm(Y) &=(-1)^{\dim Y}(c(T^*M\otimes\mathcal L)\cap s(Y,M))_{\vee\mathcal L}\\
&=(-1)^{\dim M-\dim Y} c(TM)\cap\left(s(Y,M)^\vee \otimes \mathcal
L\right)\quad.
\endaligned$$
Therefore
$$s(Y,M)^\vee\otimes\mathcal L=(1-\chi)c(N_{Y'}M)^{-1}\cap [Y']\quad.$$
Now arguing as in the proof of Lemma~\ref{lemma4}:
$$\aligned
s(Y,X)^\vee \otimes\mathcal L &=\frac{\text{Eu}-\chi}{1-\chi}\,
\frac{1+\frac{\text{Eu}-1}{\text{Eu}-\chi}X}{1+X}\,\left( s(Y,M)^\vee
\otimes \mathcal L\right)\\
&=\frac{(\text{Eu}-\chi)+(\text{Eu}-1)X}{1+X} \,c(N_{Y'}M)^{-1}\cap [Y']
\endaligned$$
and the statement follows from the expression for $\sma(X,M)$ obtained
in the proof of Theorem~\ref{main}.
\end{proof}

\begin{example}\label{example4.1}
{\em The tangent developable surface of the twisted
cubic in $\Pbb^3$.\/} For a concrete example, let $X$ be the surface
obtained as the union of all tangent lines to a fixed twisted cubic
curve in $\Pbb^3$. It is a standard but pleasant exercise to check that

---$X$ is a surface of degree~4;

---its singular locus $Y'$ is the twisted cubic;

---the polar classes of $X$ are: $[P_0]=X$; $[P_1]=3[\Pbb^1]$;
$[P_2]=0$.

The numerical invariants of $X$ are clearly constant along the twisted
cubic: indeed, $X$ is invariant under an action of $\PGL(2)$ that is
transitive along the singular locus. 
Hence $X$ is a nice hypersurface of $M=\Pbb^3$.

Note that in this example $\text{Eu}$ equals the multiplicity of $X$
along $Y'$ (for example by the fundamental formula in
\cite{8}, \S5); however, this multiplicity is not available
without a local study of $X$. Also, we do not know of any direct way
to compute $\chi$ that does not involve a deeper local study of $X$.

The global information listed above suffices however to determine both
$\text{Eu}$ and $\chi$ in this example, by means of 
Proposition~\ref{proposition7}.
Indeed, we have
$$[P]=\frac{[X]}{1+H}-\frac{3[\Pbb^1]}{(1+H)^2}=[X]-7[\Pbb^1]+10[\Pbb^2]$$
(where $H$ denotes the hyperplane class), and $c(T\Pbb^3)\cap[P]=\cma(X,M)$
by Theorem~\ref{theorem5} (a). Hence
$$(1+X)\left(\cma(X)-\cfu(X)\right)=c(T\Pbb^3)\cap 
\left((1+X)[P]-[X]\right)=9[\Pbb^1]+18[\Pbb^0]\quad.$$
Applying Proposition~\ref{proposition7} gives then
$$9[\Pbb^1]+18[\Pbb^0]=\left((\text{Eu}-\chi) + (\text{Eu}-1) 
X\right)\cdot (3[\Pbb^1]+2[\Pbb^0])\quad,$$
from which $\text{Eu}-\chi=3$ and $7\,\text{Eu}-\chi=15$. Therefore
$$\text{Eu}=2\quad,\quad\chi=-1\quad.$$
Hence $\rho=\frac 13$ here, and by Corollary~\ref{polar}
$$\csm(X)=c(T\Pbb^3)\cap \frac{\frac 13[X]+\frac 23[P]}{1+\frac 13 X} =
[X]+6[\Pbb^1]+4[\Pbb^0]\quad.$$
\end{example}

We end with perhaps the simplest example in which our formula does
{\em not\/} apply.
\begin{example}\label{example4.2} 
Consider a reduced plane curve of degree $d\ge 3$ with
exactly one node, and let $X$ be the cone in $\Pbb^3$ over this curve.
The singular locus $Y'$ of $X$ is then a line $L$, but the singularity
scheme $Y$ is `fatter' at the vertex of the cone. The invariants
considered here detect this feature of $Y$: it is not hard to check
that $(\chi,\text{Eu})=(0,2)$ at all points of $L$ but the vertex,
while $(\chi,\text{Eu})=(d(d-1)(d-2),2+2d-d^2)$ at the vertex; in
particular, these numbers are not constant along $Y'$, so $X$ is not
`nice'.

For this example we have $[P_0]=[X]$, $[P_1]=d^2-d-2$ lines through
the vertex, and $[P_2]=0$. The push-forward to $[\Pbb^3]$ of the class
$c_\alpha(X)$ defined in \S\ref{intro} is
$$d[\Pbb^2]+\left(2+4d-d^2-2\alpha\right)[\Pbb^1] + \left(4+5d-2d^2+(-4-d
-2d^2 +d^3)\alpha+2d\alpha^2\right)[\Pbb^0].$$
It is immediate to check that this expression does {\em not\/} equal
the push forward of $\csm(X)$:
$$d[\Pbb^2]+(1+4d-d^2)[\Pbb^1]+(2+3d-d^2)[\Pbb^0]$$
for any value of $\alpha$. Note however that the value $\alpha=\frac 12$ 
corresponding to $(\chi,\text{Eu})=(0,2)$ at the general point on $L$ {\em 
does\/} yield the correct term in codimension~1: indeed, the invariants jump 
on a locus of codimension~2 in $X$, so (as observed in the introduction) the 
formula in Theorem~\ref{main} is correct for all terms of lower codimension 
in~$X$.
\end{example}

\bibliographystyle{plain}
\bibliography{interpbiblio}

\end{document}